\documentclass{amsart}
\usepackage[small, nohug]{diagrams}
\usepackage[mathscr]{eucal}
\usepackage{amsfonts}
\usepackage{amsmath}
\usepackage{amsthm}
\usepackage{amssymb}
\usepackage{latexsym}
\newfont{\msam}{msam10}

\newtheorem{theorem}[]{Theorem}
\newtheorem{proposition}[]{Proposition}

\newtheorem{lemma}[]{Lemma}
\theoremstyle{definition}
\newtheorem{definition}[]{Definition}
\newtheorem{remark}[]{Remark}
\newtheorem{example}[]{Example}

\def\remark{\noindent\textbf{Remark.}}

\let\nc\newcommand

\nc{\la}{\label}

\def\bthm{\begin{theorem}}
\def\ethm{\end{theorem}}
\def\blemma{\begin{lemma}}
\def\elemma{\end{lemma}}
\def\bproof{\begin{proof}}
\def\eproof{\end{proof}}
\def\bprop{\begin{proposition}}
\def\eprop{\end{proposition}}

\def\Gr{\mbox{\rm{Gr}}^{\mbox{\scriptsize{\rm{ad}}}}}

\def\F{\mathcal F}
\def\vp{\varphi}
\def\bvp{\boldsymbol{\psi}}
\def\bP{\boldsymbol{P}}
\def\ttau{\boldsymbol{\tau}}
\def\DC{\mathscr{D}}
\def\P{\mathcal P}
\def\Z{\mathbb{Z}}
\def\F{\mathcal{F}}
\def\R{\mathscr{I}}
\def\I{\mathcal{I}}
\def\J{\mathcal{J}}
\def\V{\mathcal{Z}}
\def\O{\mathcal{O}}
\def\A{\mathbb{A}}

\def\L{\mathcal{L}}
\def\N{\mathbb{N}}

\def\D{\mathcal{D}}
\def\DO{\overline{\mathcal D}}

\def\RR{R}
\def\MO{\overline{M}}

\def\tA{A}
\def\tsigma{\tilde{\sigma}}

\def\AA{\O^{\,\otimes 2}}
\def\AAA{A^{\otimes 2}}
\def\BB{B^{\otimes 2}}
\def\DA{\mbox{\rm DR\!}^1(\O)}
\def\DB{\mbox{\rm DR\!}^1(B)}

\def\eB{B^{\mbox{\scriptsize{\rm{e}}}}}

\def\ttheta{\tilde{\theta}}

\def\btheta{\boldsymbol{\theta}}

\def\v{\mbox{\rm v}}
\def\hv{\hat{\mbox{\rm v}}}
\def\hhv{\bar{\mbox{\rm v}}}
\def\w{\mbox{\rm w}}
\def\hw{\hat{\mbox{\rm w}}}
\def\hhw{\bar{\mbox{\rm w}}}
\def\ha{\hat{a}}
\def\hha{\bar{a}}

\def\hd{\hat{d}}
\def\hhd{\bar{d}}
\def\hdel{\hat{\Delta}}
\def\sX{\bar{X}}
\def\sY{\bar{Y}}
\def\sZ{\bar{Z}}
\def\sdel{\bar{\Delta}}
\def\sw{\bar{\w}}
\def\sv{\bar{\v}}

\def\m{\mathfrak{m}}
\def\c{\mathbb{C}}
\def\CC{\mathcal{C}}

\def\bV{\boldsymbol{V}}
\nc{\Hom}{{\rm{Hom}}}
\nc{\Ext}{{\rm{Ext}}}
\nc{\HOM}{\underline{\rm{Hom}}}
\nc{\EXT}{\underline{\rm{Ext}}}
\nc{\TOR}{\underline{\rm{Tor}}}
\nc{\End}{{\rm{End}}}
\nc{\GL}{{\rm{GL}}}
\nc{\Rep}{{\rm{Rep}}}
\nc{\ad}{{\rm{ad}}}
\nc{\dlim}{\varinjlim}

\newcommand{\Mod}{{\tt{Mod}}}
\newcommand{\Mat}{{\tt{Mat}}}

\newcommand{\HH}{{\rm{HH}}}

\newcommand{\Tor}{{\rm{Tor}}}
\newcommand{\pr}{{\tt can\,}\,}

\newcommand{\Pic}{{\rm{Pic}}}
\newcommand{\Aut}{{\rm{Aut}}}
\newcommand{\Out}{{\rm{Out}}}
\newcommand{\id}{{\rm{Id}}}
\newcommand{\Der}{{\rm{Der}}}
\newcommand{\Div}{{\rm{Div}}}
\newcommand{\rk}{{\rm{rk}}}
\newcommand{\Tr}{{\rm{Tr}}}

\newcommand{\Ker}{{\rm{Ker}}}

\newcommand{\bmu}{{\boldsymbol{\mu}}}
\newcommand{\bn}{{\boldsymbol{n}}}

\newcommand{\be}{{\boldsymbol{\rm e}}}

\newcommand{\blambda}{{\boldsymbol{\lambda}}}

\newcommand{\Vi}{V_{\infty}}

\newcommand{\ei}{e_{\infty}}
\newcommand{\lambdai}{\lambda_{\infty}}

\newcommand{\into}{\,\,\hookrightarrow\,\,}

\newcommand{\onto}{\,\,\twoheadrightarrow\,\,}
\begin{document}
\title[]{Calogero-Moser Spaces over Algebraic Curves}
%
%
\author{Yuri Berest}
\address{Department of Mathematics, Cornell University, Ithaca,
NY 14853-4201, USA}
\email{berest@math.cornell.edu}
%
%
\begin{abstract}
In these notes, we give a survey of the main results of \cite{BC}
and \cite{BW}. Our aim is to generalize the geometric
classification of (one-sided) ideals of the first Weyl algebra $
A_1(\c) $ (see \cite{BW1, BW2}) to the ring $ \D(X) $ of
differential operators on an arbitrary complex smooth affine curve
$X$. We approach this problem in two steps: first, we classify the
ideals of $ \D(X) $ up to stable isomorphism, in terms of the
Picard group of $ X $; then, we refine this classification by
describing each stable isomorphism class as a disjoint union of
(certain quotients of) generalized Calogero-Moser spaces $
\CC_n(X, \I) \,$. The latter are defined as representation
varieties of deformed preprojective algebras over a one-point
extension of the ring of regular functions on $ X $ by the line
bundle $ \I $. As in the classical case, $ \CC_n(X, \I) $ turn out
to be smooth irreducible varieties of dimension $2n$.
\end{abstract}
\maketitle
\section{Introduction}
\la{intro} In recent years, there have been a number of
interesting proposals in the area of {\it smooth} noncommutative
algebraic geometry (see \cite{CQ}, \cite{KR}, \cite{G},
\cite{LeB}, \cite{vdB}, \cite{CEG}). The algebras studied in this
area are called quasi-free or formally smooth as they appear to be
`smooth' objects in the category of associative algebras (in the
same sense as the rings of functions on nonsingular affine
varieties in the category of {\it commutative} algebras). Over the
complex numbers, quasi-free algebras can be characterized
cohomologically as the ones having dimension $ \le 1 $ with
respect to Hochschild cohomology. This characterization shows that
`quasi-freeness' is a very restrictive property. Apart from
semi-simple algebras, there are basically two sources of examples:
the path algebras of quivers and the (commutative) rings of
functions on smooth affine curves\footnote{There are also a few
natural constructions, which can be used to produce new quasi-free
algebras from the old ones. For example, the class of quasi-free
algebras is closed under products and coproducts in the category
of associative algebras as well as (universal) localizations, see
\cite{CQ}.}. Most developments in the area follow a familiar
pattern of noncommutative geometry: translating geometric concepts
and intuition into algebraic language and extending these to
arbitrary quasi-free algebras. In this way, one gets
noncommutative analogues of various results of smooth algebraic
geometry.

On the other hand, quivers bring in a rich source of ideas and
constructions originating from representation theory of
finite-dimensional algebras. Reversing the logic (and alienating,
perhaps, some classically educated geometers), one may try to
apply these to {\it commutative} quasi-free algebras, i.~e. to the
ordinary curves viewed as objects of noncommutative geometry. It
is this last approach that we adopt to define the Calogero-Moser
spaces.

One construction, which plays a fundamental role in noncommutative
geometry, is that of a representation variety of an algebra: it
generalizes the variety of representations of a quiver. Another is
that of a deformed preprojective algebra \cite{CB}\,: it
generalizes the classical preprojective algebras associated to
graphs (see \cite{GP}). The third, perhaps less known, is a
one-point extension of an algebra: this abstracts the idea of
`framing' a quiver (by adding to it a distinguished new vertex `$
\infty $' and arrows from $\, \infty \,$, see, e.g., \cite{R}).
These three constructions are key ingredients of our definition of
Calogero-Moser spaces, and we will review them in some detail in
Sections~\ref{CM}.

In Section 2, after some preliminaries, we explain our
classification of ideals of $ \D(X) $ up to isomorphism in $
K_0(\D) $, and relate this to an earlier work of Cannings and
Holland \cite{CH}. The main results of this section
(Theorem~\ref{T4} and Proposition~\ref{T5}) are proved in
\cite{BW}.

In Section~\ref{CM}, we present our definition of the
Calogero-Moser spaces $ \CC_n(X, \I) $ for an arbitrary curve $X$
and a line bundle $ \I $ over $X$. We begin with the simplest
example: $ X = \A^1 $, in which case we observe that $ \CC_n(X,
\I) $ coincide with the classical Calogero-Moser spaces $ \CC_n $
(as defined in \cite{Wi}). Apart from \cite{BC1}, this observation
was the starting point for our work. The main result of this
section (Theorem~\ref{T6}) is a generalization of a well-known
theorem of Wilson \cite{Wi} on irreducibility of the
Calogero-Moser spaces. We note that the spaces $\, \CC_n(X, \I)
\,$ behave functorially with respect to $ \I $, so the quotients
$\,\overline{\CC}_n(X,\I) := \CC_n(X, \I)/\Aut_X(\I) \,$ depend
only on the class of $ \I $ in $ \Pic(X) $. We conclude
Section~\ref{CM} with an explicit description of $ \CC_n(X, \I) $,
in terms of matrices satisfying a `rank one' condition, and
illustrate our theory with a broad class of examples, including a
general plane curve.

Finally, in Section~\ref{CMMap}, we construct a natural action of
the Picard group $ \Pic(\D) $ on the reduced Calogero-Moser spaces
$\, \overline{\CC}_n(X, \I) \,$ and state our main result
(Theorem~\ref{Tmain}). This theorem provides a classification of
left ideals of $ \D $ in terms of  $ \overline{\CC}_n(X, \I) $,
which is, like in the Weyl algebra case, equivariant under the
action of $ \Pic(\D) $. The classifying map $ \omega $ from $
\overline{\CC}_n(X, \I) $ to the space of ideals $ \R(\D) $ is
induced by a certain functor from the representation category of a
deformed preprojective algebra to the category of $ \D$-modules;
in the special case when $ X = \A^1 $, this map agrees with the
Calogero-Moser map constructed in \cite{BW1, BW2}.

There still remain many questions. First of all, in the existing
literature, there are (at least) two other definitions of
Calogero-Moser spaces associated to curves. The first one, due to
Etingof (see \cite{E}, Example~2.19), is given in terms of
generalized Cherednik algebras (in the style of \cite{EG}). The
second, due to Ginzburg, employs the classical Hamiltonian
reduction (see \cite{BN}, Definition~1.2). It is more or less
clear that all three definitions should agree with each other, but
it is not clear whether there exist {\it canonical} isomorphisms
between them.

Next, there is an alternative description of torsion-free
$\D$-modules on curves, using a noncommutative version of
Beilinson's resolution (see \cite{BN}). Despite the fact that one
of the starting points for \cite{BN} was to extend \cite{BC1} to
general curves (which was also the starting point for the present
work), the precise relation between the two approaches is not
clear to us at the moment. It seems that the methods of \cite{BN}
are suitable for projective curves, while in the affine case, lead
to a much more complicated classification of ideals than in
\cite{BC} (for example, no explicit map similar to our $ \omega $
appears in this classification). Comparing the two approaches is
still an interesting problem, which we plan to discuss elsewhere.

We should also mention some generalizations. Many results of
\cite{BC} can be extended to an arbitrary (formally) smooth
algebra, so it is natural to ask whether there is a general
principle in noncommutative geometry behind our approach. On the
other hand, it might be interesting to understand the results of
\cite{BW} and \cite{BC} in purely geometric terms, with a view of
extending them to complete and analytic curves.

In the end, I would like to thank W.~Crawley-Boevey, P.~Etingof,
V.~Ginzburg, I.~Gordon, R.~Rouquier, G.~Segal, and especially my
coauthors O.~Chalykh and G.~Wilson for interesting questions and
comments. This paper evolved from notes of my talk at the
conference on Cherednik algebras in June 2007. I would like to
thank the organizers of this conference, in particular
Iain~Gordon, for inviting me to Edinburgh and giving an
opportunity to speak. This work was partially supported by NSF
grant DMS-0407502 and a LMS grant for visiting scholars.

\subsection*{Notation} Throughout this paper, $ X $ will denote a
smooth affine irreducible curve over $\c$, $\, \O = \O(X) \,$ the
ring of regular functions on $ X $, and $ \D = \D(X) $ the ring of
global (algebraic) differential operators on $ X $. Unless
otherwise specified, a module over a ring $ R $ means a left
module over $R$, and $ \Mod(R) $ stands for the category of such
modules.

\section{Rings of Differential Operators on Curves}
\la{diffop}
In this section, we state our first result (Theorem~\ref{T4}),
which gives a $K$-theoretic classification of ideals of $ \D $.

\subsection{Basic properties}
Recall that $ \D $ is a filtered algebra $\, \D = \bigcup_{k
\ge 0}\, \D_k \,$, with filtration components $\,  0 \subset \D_0
\subset \ldots \subset \D_{k-1} \subset \D_{k} \subset \ldots \,$
defined inductively by
$$
\D_k := \{\, D \in \End_{\c}\,\O\ :\ [\,D,\,f\,] \in \D_{k-1}\
\mbox{for all}\ f \in \O \,\}\ .
$$
The elements of $ \D_k $ are called {\it differential operators of
order} $ \le k $. In particular, the differential operators of
order $ 0 $ are multiplication operators by regular functions on
$X$, i.e. $\,\D_0 = \O\,$, while the differential operators of
order $ \le 1 $ are  linear combinations of functions and
(algebraic) vector fields on $X$, i.e. $ \D_1 $ is spanned by $ \O
$ and the space $ \Der(\O)$ of derivations of $ \O \,$. When $ X $
is smooth and irreducible (as we always assume in this paper), $
\O $ and $ \Der(\O) $ generate $ \D $ as an algebra, and $ \D $
shares many properties with the first Weyl algebra $
A_1(\c)$, of which it is a generalization: $\,A_1(\c) \cong
\D(\c^1)$. Thus, like $ A_1(\c)$, $\, \D $ is a simple Noetherian
domain of homological dimension $ 1 $ (see \cite{Bj}, Ch.~2).
However, unlike $ A_1(\c)$, $\,\D\,$ has a nontrivial $ K$-group.

We write $\, \DO := \bigoplus_{k=0}^\infty \D_k/\D_{k-1} \,$ for
the associated graded ring of $ \D $, which is a commutative
algebra isomorphic to the ring of regular functions on the
cotangent bundle $ T^*X $ of $ X $. If $ M $ is a $\D$-module
equipped with a $\D$-module filtration $ \{ M_{k} \} $, we also
write $\, \MO := \bigoplus_{k=0}^\infty M_k/M_{k-1} $ for the
associated graded $\DO$-module. Using the standard terminology, we
say that $\, \{M_{k}\} \,$ is {\it good} if $ \MO $ is finitely
generated.

\subsection{Stable classification of ideals}
\la{SCI} Let $K_0(X)$ and $\Pic(X)$ denote the Gro\-then\-dieck
group and the Picard group of $ X $ respectively. By definition, $
K_0(X) $ is generated by the stable isomorphism classes of
(algebraic) vector bundles on $X$, while the elements of $ \Pic(X)
$ are the isomorphism classes of line bundles. As $ X $ is affine,
we may identify $ K_0(X) $ with $ K_0(\O)$, the Grothendieck group
of the ring $\O$, and $ \Pic(X) $ with $ \Pic(\O)$, the ideal
class group of $ \O $. There are two natural maps $\, \rk: K_0(X)
\to \Z \,$ and $\, \det: K_0(X) \to \Pic(X) \,$ defined by taking
the rank and the determinant of a vector bundle respectively, and
it is well-known that, in the case of smooth curves, $\, \rk
\oplus \det:\, K_0(X) \stackrel{\sim}{\to} \Z \oplus \Pic(X)\,$ is
a group isomorphism.

Now, let $ \R(\D) $ denote the set of isomorphism classes of
(nonzero) left ideals of $\D $. Unlike $ \Pic $ in the commutative
case, $ \R(\D) $ carries no natural structure of a group. However,
since $ \D $ is a hereditary domain, $ \R(\D) $ can be identified
with the space of isomorphism classes of rank $1$ projective
modules, and there is a natural map relating $ \R(\D) $ to $
\Pic(X) $:
\begin{equation}
\la{E10}
\gamma:\, \R(\D) \stackrel{\pr}{\longrightarrow}
K_0(\D) \xrightarrow{i_*^{-1}} K_0(X)
\stackrel{\det}{\longrightarrow} \Pic(X) \ .
\end{equation}
Here, $ K_0(\D) $ denotes the Grothendieck group of the ring $ \D
$, `$ \pr\!\! $' is the canonical map assigning to the isomorphism
class of an ideal of $ \D $ its stable isomorphism class in $
K_0(\D) $, and $\, i_*^{-1} \,$ is the inverse of the Quillen
isomorphism $\,\, i_*:\, K_0(X) \stackrel{\sim}{\to} K_0(\D)\,$
induced by the natural inclusion $\, i:\,\O \into \D\,$ (see
\cite{Q}, Theorem~7).

The role of the map $ \gamma $ becomes clear from the following
theorem, which is the main result of this section.
\begin{theorem}[\cite{BW}, Proposition~2.1]
\la{T4} Let $ M $ be a projective $\D$-module of rank $1$ equipped
with a good filtration. Assume that $ \MO $ is
torsion-free\footnote{This assumption is very restrictive: if we
identify $ M $ with an ideal in $ \D $, the given filtration on $
M $ coincides, up to a shift, with the induced filtration $\,\{M
\cap \D_k \} \,$ (see \cite{BC}, Lemma~5.12).}. Then

$(a)$ there is a unique (up to isomorphism) ideal $ \I_{M}
\subseteq \mathcal{O} \,$, such that $ \MO $ is isomorphic to a
sub-$\DO$-module of $  \DO \I_M $ of finite codimension (over
$\c$);

$(b)$ the class of $ \I_M $ in $ \Pic(X) $ and the codimension $\,
n := \dim_{\c}\,[\DO \I_M/\MO]\,$ are independent of the choice of
filtration on $ M $, and we have $\,\gamma [M] = [\I_M] \,$;

$(c)$ if $ M $ and $ N $ are two projective $\D$-modules of rank
$1$, then
$$
[\,M\,] = [\,N\,] \ \mbox{in}\ K_0(\D)
\quad \Longleftrightarrow \quad [\I_{M}] = [\I_{N}]
\ \mbox{in}\ \Pic(X)\ .
$$
\end{theorem}

It is easy to see that $\,\gamma[\D \I] = [\I] \,$ for any nonzero
ideal $ \I \subseteq \O $. Thus, by Theorem~\ref{T4}, the map $
\gamma $ is a fibration over $ \Pic(X) $, with fibres being
precisely the {\it stable} isomorphism classes of ideals of $ \D
$. The stably free ideals $ M $ are characterized by the property
that $ \MO $ is isomorphic to an ideal in $ \DO $ of finite
codimension.

\subsection{The Cannings-Holland correspondence} By a theory of
Cannings and Holland (see \cite{CH}), the ideals of $ \D $ can be
parametrized by primary decomposable subspaces of $ \O $. We now
describe the map \eqref{E10} in terms of these subspaces.

First, we recall that a linear subspace $\, V \subseteq \O \,$ is
called {\it primary} if it contains a power of the maximal ideal $
\m_x $ corresponding to a point $ x \in X $ (we write $ V = V_x $
in this case). Further, $\, V \subseteq \O(X) $ is called {\it
primary decomposable} if it is an intersection of primary
subspaces $\, V_x \,$, with $ V_x = \O(X) $ for almost all $ x \in
X $. By \cite{CH}, Theorem 2.4, the primary decomposition of $ V $
is uniquely determined: in fact, we have $\, V_x =
\bigcap_{k=1}^{\infty}(V + \m_x^k) \,$ for all $ x \in X $, and
moreover, $\, \dim_{\c} \O/V = \sum_{x \in X} \dim_{\c} \O/V_x
\,$.

Now, let $ M $ be a nonzero left ideal of $ \D$. Then, by
\cite{CH}, Theorem~4.12, there is a unique (up to rational
equivalence) primary decomposable subspace $\, V \subseteq \O \,$,
such that $\,M \cong \D(V, \O)\,$, where $\, \D(V, \O) \,$ is the
{\it fractional} ideal of $ \D $ consisting of all differential
operators with rational coefficients on $ X $ mapping $ V $ into
$ \O $. We write $ V_x $ for the primary
components of $V $, and $ \m_x \subset \O $ for the associated primes.
\begin{proposition}[\cite{BW}, Theorem~5.2]
\la{T5} The map $ \gamma $ sends the class of $ M $ to the class
of the ideal $\, \prod_{x \in X} \m_{x}^{d_x} \,$ in $ \Pic(X) $,
where $ d_x := \dim_{\c} \O/V_x $.
\end{proposition}

Let $ \Gr(X)$ be the {\it adelic Grassmannian} of $ X $, i.~e. the
set of equivalence classes of primary decomposable subspaces of $
\O(X)$. There is a well-defined map from $\, \Gr(X)\,$ to the
divisor class group of $ X $: it takes the class of a subspace
$\,V = \bigcap_{x \in X} V_x\,$ in $ \Gr(X) $ to the class of the
(Weil) divisor $\, d := \sum_{x \in X} d_x \cdot x \,$ in $\,
\Div(X)\,$. On the other hand, by the Cannings-Holland Theorem, we
have the bijection: $\, \Gr(X) \stackrel{\sim}{\to} \R(\D)\,$,$\
[V] \mapsto [\D(V, \O)]\,$, and, as $ X $ is smooth, the natural
isomorphism: $\, \Div(X) \stackrel{\sim}{\to} \Pic(X) \,$, $\,
[\,d\,] \mapsto [\O_X(d)] \,$. In this way, we get the diagram
\begin{equation}
\la{D5}
\begin{diagram}[small, tight]
\Gr(X)           & \rTo  & \Div(X)  \\
\dTo^{\cong}             &       & \dTo_{\cong}     \\
\R(\D)           & \rTo^{\!\!\!\! \gamma}  &  \Pic(X)
\end{diagram}
\end{equation}

\noindent Now, Proposition~\ref{T5} immediately implies that
\eqref{D5} is commutative. This gives an alternative description
of our map $\,\gamma\,$ in terms of primary decomposable
subspaces:
$$
\Gr(X) \to \Div(X)\ ,\quad [\,V\,] \mapsto
\sum\, d_x \cdot x \ ,
$$
where $ d_x $ are codimensions of the primary components of $ V $.

\section{The Calogero-Moser Spaces}
\la{CM} Theorem~\ref{T4}  shows that the ideals of $ \D(X) $ are
classified, up to stable isomorphism, by the elements of $ \Pic(X)
$. Our goal now is to refine this classification by describing the
fibres of the classifying map $\, \gamma:\, \R(\D) \to \Pic(X)
\,$. As we will see, each fibre of $ \gamma $ breaks up into a
countable union of the quotient spaces $ \overline{\CC}_n(X,\I) =
\CC_n(X,\I)/\Aut_X(\I)\,$ of smooth affine varieties $ \CC_n(X,\I)
$. The varieties $ \CC_n(X,\I) $ will be introduced as
representation varieties of deformed preprojective algebras over
the one-point extension of the ring of regular functions on $X$ by
the line bundle $\I$. In the special case when $ X $ is the affine
line, $\, \CC_n(X,\I) \,$ coincide with the ordinary
Calogero-Moser spaces \cite{Wi}, and our classification of ideals
of $ \D(X) $ agrees with the one given in \cite{BW1, BW2}.

We begin by reviewing the basic ingredients of our construction.

\subsection{Representation varieties}
\la{modvar} First, we recall the definition of representation
varieties in the form they appear in representation theory of
associative algebras (see, e.g., \cite{K}, Chap.~II, Sect.~2.7).

Let $ R $ be a finitely generated associative $\c$-algebra, $\,S
\,$, a finite-dimensional semi-simple subalgebra of $ R $, and $ V
$, a finite-dimensional $S$-module. By definition, the {\it
representation variety} $\, \Rep_{S}(R,V) \,$ of $ R $ over $ S $
parametrizes all $R$-module structures on the vector space $V$
extending the given $S$-module structure on it. The $S$-module
structure on $ V $ determines an algebra homomorphism $ S \to
\End(V) $ making $ \End(V) $ an $S$-algebra. The points of
$\,\Rep_{S}(R,V)\,$ can thus be interpreted as $S$-algebra maps
$\, R \to \End(V) $.

If $ S = \c $, we simply write $ \Rep(R, V) $ for $ \Rep_{\c}(R,
V) $. Choosing a basis in $ V $ and a presentation of $ R $, say $
R \cong \c\langle x_1, \ldots, x_m \rangle/I $, we can identify in
this case
$$
\Rep(R, V) \cong \{(X_1, \ldots, X_m) \in
\Mat(n,\c)^m:\, r(X_1, \ldots, X_m) = 0 \, ,\ \forall\,
r \in I\}\ .
$$
Thus $ \Rep(R, V) $ is an affine variety\footnote{Here, by an
affine variety we mean an affine scheme of finite type over $
\c$.}. In general, for any semi-simple $\,S \subseteq R\,$,
$\,\Rep_{S}(R,V) \,$ can be identified with a fibre of the
canonical morphism of affine varieties $\,\pi:\,\Rep(R, V) \to
\Rep(S,V)\,$, and hence it is an affine variety as well.

The group $ \Aut_S(V) $ of $S$-linear automorphisms of $ V $ acts
on $\, \Rep_{S}(R,V) \,$ in the natural way, with scalars $\,
\c^{\times} \subseteq \Aut_S(V) \,$ acting trivially. We set $\,
\GL_S(V) := \Aut_S(V)/\c^\times $. Since $ V $ is semi-simple, $\,
V \cong \bigoplus_{i} V_i^{\oplus n_i} $, with $ V_i $
non-isomorphic simple $S$-modules, and $\, \Aut_S(V) \cong
\prod_{i} \GL(n_i, \c) \,$. Thus $\, \GL_S(V) \,$ is reductive.

The orbits of $\GL_S(V) $ on $\, \Rep_{S}(R,V) \,$ are in 1-1
correspondence with isomorphism classes of $R$-modules, which are
isomorphic to $\, V \,$ as $S$-modules. The stabilizer of a point
$\, \varrho:\,R \to \End(V) \,$ in $\, \Rep_{S}(R,V) \,$ is
canonically isomorphic to $ \Aut_R(V_\varrho)/\c^\times \subseteq
\GL_S(V)$, where $ V_\varrho $ is the left $R$-module
corresponding to $ \varrho $.

Now, one can show that the closure of any orbit $ {\mathcal O}_{M}
$ contains a unique {\it closed} orbit, corresponding to a
semi-simple $R$-module with the same composition factors and
multiplicities as $ M $. Thus the space $\,
\Rep_{S}(R,V)/\!/\GL_S(V) \,$ of closed orbits in $\,
\Rep_{S}(R,V) \,$ is an affine variety, whose (closed) points are
in bijection with isomorphism classes of {\it semi-simple}
$R$-modules $M$ isomorphic to $\, V \,$ as $S$-modules.

Typically, the representation varieties of $R$ are defined over
subalgebras spanned by idempotents. For example, let  $\,
\{e_i\}_{i \in I} \,$ be a complete set of orthogonal idempotents
in $R$. Set $\, S := \bigoplus_{i \in I} \c \,e_i \subseteq R \,$.
A finite-dimensional $S$-module is then isomorphic to a direct sum
$\, \c^{\bn} := \bigoplus_{i \in I} \c^{n_i} $, each $\, e_i \,$
acting as the projection onto the $i$-th component. The
corresponding representation variety $\,\Rep_{S}(R,\c^\bn)\,$,
which we denote simply $ \Rep_S(R, \bn) $ in this case, consists
of all algebra maps $\, R \to \End(\c^{\bn})$, sending $ e_i $ to
the projection onto $ \c^{n_i} $. The group $ \GL_S(\c^\bn) $ (to
be denoted $\, \GL_S(\bn) \,$) is isomorphic to $\,\prod_{i\in I}
\GL(n_i, \c)/\c^\times \,$, with $ \c^\times $ embedded
diagonally.

\subsection{Deformed preprojective algebras}
\la{DPA} If $ A $ is an associative algebra, its tensor square $
\AAA $ over $ \c $ has two commuting bimodule structures: one is
defined by $\, a.(x \otimes y).b = ax \otimes yb \,$ and the other
by $\, a.(x \otimes y).b = xb \otimes ay\,$. We will refer to the
first structure as {\it outer} and to the second as {\it inner}.

The space $\,\Der(A, \AAA) \,$ of linear derivations $\, A \to
\AAA \,$ taken with respect to the outer bimodule structure on $
\AAA $ is naturally a bimodule with respect to the inner
structure; thus, we can form the tensor algebra $\, T_A \Der(A,
\AAA) \,$. If $A$ is unital, there is a canonical element in
$\Der(A, \AAA)$: namely the derivation $\, \Delta = \Delta_A: A
\to \AAA \,$, sending  $ x \in A $ to $ (x \otimes 1 - 1 \otimes
x) \in \AAA $. For any $\, \lambda \in A \,$, we can form then the
two-sided ideal $\,\langle \Delta - \lambda \rangle \,$ in
$\,T_A\Der(A, \AAA)\,$ and, following \cite{CB}, define the
quotient algebra
\begin{equation}
\la{E19}
\Pi^{\lambda}(A) := T_A \Der(A, \AAA)/\langle \Delta - \lambda \rangle \ .
\end{equation}
It turns out that, up to isomorphism, $\,\Pi^\lambda(A)\,$ depends
only on the class of $ \lambda $ in the Hochschild homology  $\,
\mbox{\rm HH}_0(A) := A/[A,A]\,$ (see \cite{CB}, Lemma~1.2).
Moreover, instead of elements of $ \mbox{\rm HH}_0(A) $, it is
convenient to parametrize the algebras \eqref{E19} by the elements
of $\, \c \otimes_{\Z} K_0(A) \,$, relating this last vector space
to $ \mbox{\rm HH}_0(A) $ via a Chern character map. To be
precise, let $\, \Tr_A: K_0(A) \to \mbox{\rm HH}_{0}(A) \,$ be the
map, sending the class of a projective module $ P $ to the class
of the trace of any idempotent matrix $\, e \in \Mat(n, A)\,$,
satisfying $\, P \cong e A^{\oplus n} $. By additivity, this
extends to a linear map $\, \c \otimes_{\Z} K_0(A) \to \mbox{\rm
HH}_0(A) \,$ to be denoted also $ \Tr_A $. Following \cite{CB}, we
call the elements of $\, \c \otimes_{\Z} K_0(A) \,$ {\it weights}
and define the {\it deformed preprojective algebra of weight} $\,
\blambda \in \c \otimes_{\Z} K_0(A)\,$ by
\begin{equation}
\la{E20} \Pi^{\blambda}(A) := T_A \Der(A,
\AAA)/\langle \Delta - \lambda \rangle \ ,
\end{equation}
where $\, \lambda \in A \,$ is any lifting of
$\,\Tr_A(\blambda)\,$ to $A$. Note, if $ A $ is commutative, then
$\, \mbox{\rm HH}_0(A) = A \,$, and $\, \lambda \,$ is uniquely
determined by $\,\Tr_A(\blambda)\,$.

For basic properties and examples of the algebras $
\Pi^{\blambda}(A) $, we refer the reader to \cite{CB}. Here, we
only review one important example and two theorems, which play a
role in our construction.
\begin{example}[see \cite{CB}, Lemma~4.1 and Theorem~0.3]
\la{TCB2} If $ A = \O(X) $ is the coordinate ring of a smooth
affine curve $X$, then $A$ is quasi-free. In this case, $\,
\Pi^0(A) \,$ is isomorphic to the coordinate ring $ \O(T^*X) $ of
the cotangent bundle of $ X $, and $\, \Pi^{1}(A) $ to the
filtered algebra $\D(X) $ of differential operators\footnote{We
warn the reader that, unlike $ \Pi^0(A) \cong \O(T^*X)$, the
isomorphism $ \Pi^1(A) \cong \D(X) $ is not canonical. To remedy
this, one should replace $ \D(X) $ by the ring $
\D(\Omega_X^{1/2})$ of twisted differential operators on a square
root of the canonical bundle of $X$, i.~e. on half-forms on $X$.
The existence of a canonical isomorphism $\,\Pi^1(A) \cong
\D(\Omega_X^{1/2}) \,$ was first noticed by V.~Ginzburg (see
\cite{G}, Sect. 13.4); however, the proof sketched in \cite{G} was
not quite complete. A simple argument clarifying Ginzburg's proof
can be found in the Appendix of \cite{BC} written by G.~Wilson.}.
\end{example}

In general, the $ \Pi^\blambda$-construction is not functorial in
$A$; however, it does behave well with respect to a class of
algebra maps called {\it pseudo-flat}\footnote{These are, in a
sense, `smooth' morphisms in the category of associative
algebras.}. To be precise, the pseudo-flat algebra homomorphisms
$\,\theta: B \to A \,$  are characterized by the condition:
$\,\Tor^B_1(A, A) = 0\,$, and the functoriality of $ \Pi^\blambda
$ is stated as follows.

\begin{theorem}[\cite{CB}, Theorem~0.7]
\la{TCB1} If $\, \theta: B \to A \,$ is a pseudo-flat ring
epimorphism, then, for any $\, \blambda \in \c \otimes_{\Z} K_0(B)
\,$, there is a canonical algebra map $\, \btheta:
\Pi^{\blambda}(B) \to \Pi^{\theta^{*}(\blambda)}(A)\,$, where $\,
\theta^*:\, \c \otimes_{\Z} K_0(B) \to \c \otimes_{\Z} K_0(A)\,$
is a linear map induced functorially by $\, \theta \,$. If $ B $
is quasi-free and finitely generated, then $\, \btheta \,$ is also
a pseudo-flat ring epimorphism.
\end{theorem}

Finally, the last theorem, which we want to state in this section,
provides a simple homological principle for studying
representations of $ \Pi^\blambda(A) $. It plays an important role
in \cite{BC}, underlying several proofs and constructions.

\begin{theorem}[\cite{BC}, Theorem~2.2]
\la{lift} Let $ A $ be a finitely generated quasi-free algebra,
and let $\, \varrho: \, A \to \End(V) \,$ be a representation of $
A $ on a (not necessarily finite-dimensional) vector space $ V $.
Then $ \varrho $ can be extended to a representation of $
\Pi^{\blambda}(A) $ if and only if the homology class of $\,
\varrho(\blambda) \,$ in $\, \HH_0(A,\,\End\,V) \,$ is zero. If it
exists, an extension of $\, \varrho \,$ to $\, \Pi^{\blambda}(A)\,
$ is unique if and only if $\, \HH_1(A,\,\End\,V) = 0 \,$.
\end{theorem}
\remark\  To the best of our knowledge, Theorem~\ref{lift} has not
appeared in the literature in this form and generality. However,
in the special case when $ A = \c Q $ is the path algebra of a
quiver and $ V $ is a {\it finite-dimensional} representation of $
A $, this result is equivalent to \cite{CB2}, Theorem~3.3 (see
\cite{BC}, Proposition~2.1).

\subsection{One-point extensions}
\la{1}
If $ A $ is a unital associative algebra, and $ I $ a left
module over $ A $, we define the {\it one-point extension} of $ A
$ by $ I $ to be the ring of triangular matrices
\begin{equation}
\la{E16} A[I] := \left(
\begin{array}{cc}
A & I \\
0 & \c
\end{array}
\right)
\end{equation}
with matrix addition and multiplication induced from the module
structure of $ I $. Clearly, $ A[I] $ is a unital associative
algebra, with identity element being the identity matrix. There
are two distinguished idempotents in $ A[I]\, $: namely
\begin{equation}
\la{idem01} e := \left(
\begin{array}{cc}
 1  & 0 \\
 0  & 0
\end{array}
\right) \quad \mbox{and} \quad \ei := \left(
\begin{array}{cc}
 0 & 0 \\
 0 & 1
\end{array}
\right) \ .
\end{equation}
If $ A $ is indecomposable (e.g., $A$ is a commutative integral
domain), then \eqref{idem01} form a complete set of primitive
orthogonal idempotents in $ A[I] $.

A module over $ A[I] $ can be identified with a triple $\, \bV =
(V,\, \Vi, \, \varphi) \,$, where $ V $ is an $A$-module, $ \Vi $
is a $\c$-vector space and $ \varphi: \, I \otimes \Vi \to V $ is
an $A$-module map. Using the standard matrix notation, we will
write the elements of $ \bV $ as column vectors $ (v, w)^T $ with
$ v \in V $ and $ w \in \Vi $; the action of $ A[I] $ is then
given by
$$
\left(
\begin{array}{cc}
a & b \\
0 & c
\end{array}
\right) \left(
\begin{array}{c}
v \\
w
\end{array} \right)
= \left(
\begin{array}{c}
a.v + \varphi(b \otimes w) \\
c w
\end{array} \right)\ .
$$
If $ \bV $ is finite-dimensional, with $\,\dim_{\c} V = n\,$ and
$\,\dim_{\c} \Vi = n_{\infty}\,$, we call $\,\bn = (n,\,
n_{\infty}) \,$ the {\it dimension vector} of $ \bV $.

The next lemma collects some basic properties of one-point
extensions, which we will need for our construction.
\blemma \la{LL6} $(1)\,$ $ A[I] $ is canonically isomorphic to the
tensor algebra $\, T_{A \times \c}(I) \,$.

$(2)\,$  If $ A $ is quasi-free and $ I $ is f.~g. projective,
then $ A[I] $ is quasi-free.

$(3)\,$ $\, I \mapsto A[I] \,$ is a functor from $ \Mod(A) $ to
the category of associative algebras.

$(4)\,$ The natural projection $\,\theta:\,A[I] \to A \,$ is a
flat ring epimorphism.

$(5)\,$ There is an isomorphism of abelian groups $ K_0(A[I])
\cong K_0(A) \oplus \Z \,$. \elemma

\noindent For the proof of Lemma~\ref{LL6}, we refer the reader to
\cite{BC}, Section~2.2.

\subsection{The definition of Calogero-Moser spaces}
\la{GCMS} We can now put pieces together and introduce our
generalization of the Calogero-Moser varieties for an arbitrary
smooth affine curve $ X $. We begin with the simplest example: $\,
X = \A^1 \,$, which will provide motivation for our general
construction.

\begin{example}
\la{Ex1} If $ X = \A^1 $, any line bundle $ \I $ on $ X $ is
trivial. Choosing a coordinate on $ X $ and a trivialization of $
\I $, we identify $\, \O \cong \c[x]\,$ and $\, \I \cong
\c[x]\,$ as an $\O$-module. The one-point extension of $ \O $ by $
\I $ is then isomorphic to the matrix algebra:
$$
\O[\I] \cong \left(
\begin{array}{cc}
\c[x] & \c[x]\\*[1.1ex]
0    & \c
\end{array}
\right)\ ,
$$
which is, in turn, isomorphic to the path algebra $ \c Q $ of the
quiver $ Q $ consisting of two vertices $\,\{0, \,\infty\}\,$ and
two arrows $\, X:\, 0 \to 0 \,$ and $\, \v:\, \infty \to 0 $. In
fact, the map sending the vertices $0$ and $\infty $ to the
idempotents $ e $ and $ \ei $ in $ \O[\I ] $, see \eqref{idem01},
and
$$
X \mapsto \left(
\begin{array}{cc}
x & 0\\
0 & 0
\end{array}
\right)\quad , \qquad
 \v \mapsto \left(
\begin{array}{cc}
0 & 1\\
0 & 0
\end{array}
\right)\ ,
$$
extends to an algebra isomorphism $\,\c Q \stackrel{\sim}{\to}
\O[\I]\,$.

Now, let $ \bar{Q} $ be the double quiver of $ Q $ obtained by
adding the reverse arrows $ Y := X^* $ and $ \w := \v^* $ to the
corresponding arrows of $Q$. Then, for any $\,\blambda =
\lambda\,e + \lambdai \ei \,$, with $ (\lambda, \,\lambdai) \in
\c^2 $, the algebra $\, \Pi^{\blambda}(Q) \,$ is isomorphic to the
quotient of $ \c\bar{Q} $ modulo the relation $\, [X,Y] + [\v, \w]
= \blambda \,$ (see \cite{CB}, Theorem~3.1). The ideal generated
by this last relation is the same as the ideal generated by the
elements $\, [X,Y] + \v\,\w - \lambda \,e \,$ and $\, \w\,\v
+\lambdai \ei \,$. Thus, the $ \Pi^{\blambda}(Q)$-modules can be
identified with representations $ \bV = V \oplus \Vi $ of $
\bar{Q} $, in which linear maps $ \,\sX, \sY \in \Hom(V,\, V),
\,\sv \in \Hom(\Vi,\, V),\, \sw \in \Hom(V,\, \Vi) \,$, given by
the action of $\, X, Y, \v, \w \,$, satisfy
\begin{equation}
\la{E26} [\sX,\,\sY] + \sv \,\sw  = \lambda \, \id_{V}\quad
\mbox{and} \quad \sw\,\sv = -\lambdai \, \id_{\Vi} \ .
\end{equation}
Now, taking $ \blambda = (1, -n) $, it is easy to see that all
representations of $ \Pi^\blambda(Q) $ of dimension vector $ \bn =
(n,\,1) $ are simple, and the corresponding representation
varieties coincide (in this special case) with the classical
Calogero-Moser spaces $ \CC_n $. This coincidence was first
noticed by W.~Crawley-Boevey (see \cite{CB1}, remark on p.~45).
For explanations and a detailed discussion of this example in
relation to the Weyl algebra we refer the reader to \cite{BCE}.
\end{example}

Now, let $ X $ be an arbitrary curve. As in the above example, we
fix a line bundle $ \I $ on $ X $ and set $ B := \O[\I] $. Note
that, by Lemma~\ref{LL6}$(3)$, $\, B \,$ depends (up to
isomorphism) only on the class of $ \I $ in $ \Pic(X) $. More
precisely, we have
\begin{lemma}[\cite{BC}, Lemma~3.1]
\la{L10} For two line bundles $ \I $ and $ \J $ on $X$, the
algebras $ \O[\I] $ and $ \O[\J] $ are

$(a)$ Morita equivalent;

$(b)$ isomorphic if and only if $ \J \cong \I^\tau $ for some $
\tau \in \Aut(X) $, where $ \I^\tau := \tau^*\I$.
\end{lemma}

To define the deformed preprojective algebras over $ B $ we need
to compute the Chern character $\,\Tr_B : K_0(B) \to \HH_0(B)\,$.
We recall that $\,\Tr_{*} :\, K_0 \to \HH_0 \,$ is a natural
transformation of functors on the category of associative
algebras: corresponding to the algebra map $\,\ttheta: B \to \O
\times \c\,$, we have thus a commutative diagram
\begin{equation}
\la{D6}
\begin{diagram}[small, tight]
K_0(B)             & \rTo^{\Tr_B}     & \HH_0(B) \\
\dTo^{}&                  & \dTo_{} \\
K_0(\O \times \c)   & \rTo^{\Tr_{\O \times \c}\ } &
\HH_0(\O \times \c) \\
\end{diagram}
\end{equation}
The two vertical maps in \eqref{D6} are isomorphisms: the first
one is given by Lemma~\ref{LL6}$(5)$, while the second has the
obvious inverse (induced by the embedding $\, \O \times \c \into B
\,$ via diagonal matrices). We will use these isomorphisms to
identify $\,\HH_0(B) \cong \HH_0(\O \times \c) \cong \O \oplus \c
\subset B \,$ and
\begin{equation}
\la{K}
K_0(B) \cong K_0(\O \times \c) \cong \Z \oplus \Z \oplus \Pic(X)\ ,
\end{equation}

Now, for any commutative algebra (e.g., $ \O \times \c $), the
Chern character map factors through the rank. Hence, with above
identifications, $\,\Tr_{B}\,$ is completely determined by its
values on the first two summands in \eqref{K}, while vanishing on
the last. Since $\,\Tr_B[(1,0)] = e \,$ and $\,\Tr_B[(0,1)] = \ei
\,$, the linear map $\, \Tr_B: \c \otimes_{\Z} K_0(B) \to \HH_0(B)
\,$ takes its values in the two-dimensional subspace $ S $ of $ B
$ spanned by the idempotents $ e $ and $ \ei $. Identifying $ S $
with $ \c^2$, we may regard the vectors $\, \blambda :=
(\lambda,\, \lambdai) = \lambda e + \lambdai \ei \in S \,$ as
weights for the family of deformed preprojective algebras
associated to $ B \,$:
\begin{equation}
\la{dpa1}
\Pi^{\blambda}(B) = T_{B}\Der(B, \BB)/ \langle \Delta_B
- \blambda \rangle\ .
\end{equation}

Since $ \O $ is an integral domain, $\,\{e,\, \ei\}\,$ is a
complete set of primitive orthogonal idempotents in $
\Pi^{\blambda}(B) $, and $\, S = \c \,e \oplus \c\, \ei \,$ is the
associated semi-simple subalgebra of $ \Pi^{\blambda}(B) $.

For each $\, \bn = (n,\, n_\infty) \in \N^2 $, we now form the
representation variety $\, \Rep_{S}(\Pi^{\blambda}(B), \,\bn) \,$
over $ S $ and define
\begin{equation}
\la{E25}
\CC_{\bn, \blambda}(X, \I) := \Rep_{S}(\Pi^{\blambda}(B),
\bn)/\!/ \GL_S(\bn)\ .
\end{equation}
As explained in Section~\ref{modvar}, $\, \CC_{\bn, \blambda}(X,
\I) \,$ is an affine scheme, whose (closed) points are in
bijection with isomorphism classes of semi-simple
$\Pi^{\blambda}(B)$-modules of dimension vector $ \bn $.

Now, by Lemma~\ref{LL6}$(3)$, every automorphism of $ \I $ induces
an $S$-algebra automorphism of $\,\Pi^{\blambda}(B)\,$  and hence
an automorphism of the representation variety $\,\CC_{\bn,
\blambda}(X, \I)\,$. We let $ \overline{\CC}_{\bn, \blambda}(X,
\I) $ denote the corresponding quotient space:
\begin{equation}
\la{quots}
\overline{\CC}_{\bn, \blambda}(X, \I) :=
\CC_{\bn, \blambda}(X, \I)/\Aut_X(\I)\ .
\end{equation}
By definition, $\, \overline{\CC}_{\bn, \blambda}(X, \I) \,$
depends only on the class of $ \I $ in $ \Pic(X) $.

More generally, from Lemma~\ref{L10}$(1)$ and the fact that
$\,\Pi^\blambda $ behaves naturally under Morita equivalence (see
\cite{CB}, Corollary~5.5), it follows that
$$
\CC_{\bn, \blambda}(X, \I) \cong\CC_{\bn, \blambda}(X, \J)
$$
for {\it any} line bundles $ \I $ and $ \J $; however, there is no
natural choice for such an isomorphism.

Motivated by the above example, we will be interested in
representations of $ \Pi^{\blambda}(B)$ of dimension $ \bn = (n,1)
$. Using Theorem~\ref{lift}, it is not difficult to prove that
such representations may exist only if $\, \blambda = 0 \,$ or $\,
\blambda = (\lambda, -n\lambda) $, with $ \lambda \not=0 $. In
this last case, the algebras $\Pi^{\blambda}(B) $ are all
isomorphic to each other, so without loss of generality we may
take $\lambda=1$.

\begin{proposition}[\cite{BC}, Proposition~3.2]
\la{L11} Let $ \blambda = (1, -n) $ and $ \bn = (n,1) $ with $ n
\in \N $. Then, for any $\,\I\,$, $\Pi^{\blambda}(B) $ has
representations of dimension vector $ \bn $, and every such
representation is simple.
\end{proposition}

We are now in position to state the main definition and the main
theorem of this section.

\vspace{1ex}

\begin{definition}
\la{CMv} The variety $ \CC_{\bn, \blambda}(X, \I) $ with
$\,\blambda = (1,-n) $ and $ \bn = (n,1) $ will be denoted $
\CC_n(X,\I) $ and called the {\it $n$-th Calogero-Moser space of
type} $ (X, \I) $. The corresponding quotient \eqref{quots} will
be denoted $\,\overline{\CC}_{n}(X, \I)\,$ and called the $n$-th
{\it reduced} Calogero-Moser space.
\end{definition}

\vspace{1ex}

In view of Proposition~\ref{L11}, the varieties $\, \CC_n(X, \I) $
parametrize the isomorphism classes of simple $
\Pi^{\blambda}(B)$-modules of dimension $ \bn = (n,1) $; they are
non-empty for any $\, [\I\,] \in \Pic(X) \,$ and $ n \geq 0 $. In
the special case when $ X $ is the affine line, $\, \CC_n(X,\I) $
coincide with the ordinary Calogero-Moser spaces $ \CC_n $ as
defined in \cite{Wi} (see Example~\ref{Ex1}). Now, one of the main
results of \cite{Wi} says that each $\CC_{n} $ is a smooth affine
irreducible variety of dimension $ 2n $. The following theorem
shows that this is true in general.
\begin{theorem}[\cite{BC}, Theorem~3.2]
\la{T6} For each $ n \geq 0 $ and $ [\I\,]\in \Pic(X) $, $\,
\CC_{n}(X,\I) \,$ is a smooth irreducible affine variety of
dimension $ 2n $.
\end{theorem}

We close this section by describing {\it generic} points of the
varieties $ \CC_n(X, \I) $ in geometric terms. First of all, using
Theorem~\ref{lift}, it is not difficult to show that any $
\Pi^\blambda(B) $-module of dimension vector $ \bn = (n,1) $
restricts to an indecomposable $B$-module, and conversely, every
indecomposable $B$-module of dimension vector $ \bn = (n,1) $
extends to a $ \Pi^\blambda(B) $-module. The generic points of $
\CC_n(X, \I) $ correspond under this restriction/extension to the
$B$-modules $ \bV $ with $\, \End_B(\bV) \cong \c \,$. Now, as
explained in Section~\ref{1}, a $B$-module structure on $ \bV = V
\oplus \Vi $ is determined by an $\O$-module homomorphism $\,
\varphi:\, \I \otimes \Vi \to V \,$, and if $\, \dim_{\c} \Vi = 1
$, it is easy to see that $\, \End_B(\bV) \cong \c \,$ is
equivalent to $ \varphi $ being surjective. Thus, for constructing
generic points of $ \CC_n(X,\I)$, it suffices to construct a
torsion $ \O$-module $V$ on $ X $ of length $ n $ together with a
{\it surjective} $\O$-module map $ \varphi: \I \to V $.
Geometrically, this can be done as follows.

Identify $ \I $ with an ideal in $ \O $ and fix $\, n \,$ {\it
distinct} points $\,p_1,\,p_2,\,\ldots,\,p_n \,$ on $ X $ outside
the zero locus of $\,\I\,$. Let $\, V := \O/\J\,$, where $\,\J\,$
is the product of the maximal ideals $ \m_i \subset \O $
corresponding to $p_i$'s. Clearly, $\, \O/\J \cong \bigoplus_{i=1}^n
\O/\m_i \, $ and $ \dim_{\c} V = n $. Now, since $ \O $ is a
Dedekind domain and $\, \I \not\subset \m_i \,$ for any $ i =
1,\,2,\,\ldots, \,n\,$, we have $\, (\O/\J) \otimes_\O (O/\I)
\cong \bigoplus_{i=1}^n (\O/\m_i) \otimes_\O (\O/\I)= 0 \,$ and $\,
\Tor_1(\O/\J, \O/\I) \cong (\I \cap \J)/\I\J = 0 \,$, so the
canonical map $\, V \otimes_\O \I  \to V $ is an isomorphism. On
the other hand, as $ V $ is a cyclic $ \O$-module, $ \I$ surjects
naturally onto $\, V \otimes_\O \I \,$. Combining $\, \I \onto V
\otimes_\O \I \stackrel{\sim}{\to} V \,$, we get the required
homomorphism $ \varphi $.

\subsection{The structure of $ \Pi^\blambda(B) $}
\la{strP} One advantage of defining the Calogero-Moser spaces as
representation varieties is that they can be described explicitly,
like in the classical case, in terms of matrices satisfying the
`rank-one condition' (see \eqref{r0} below). For this, it suffices
to find a suitable presentation of the algebras $ \Pi^\blambda(B)
$ in terms of generators and relations.

Recall that, following \cite{CB}, we defined these algebras by
$$
\Pi^{\blambda}(B) = T_{B}\Der(B, \BB)/ \langle \Delta_B - \blambda
\rangle\ ,
$$
where $\,\Delta_B \,$ is the distinguished derivation in
$\,\Der(B, \BB) \,$ mapping $\, x \mapsto x \otimes 1 - 1 \otimes
x \,$. Now, $\,\Der(B, \BB)\,$ contains a canonical sub-bimodule
$\, \Der_S(B, \BB)\,$, consisting of $S$-linear derivations. We
write $\, \Delta_{B,S}: \,B \to B \otimes B \,$ for the inner
derivation $\, x \mapsto \ad_{\be}(x) \,$, with $\, \be := e
\otimes e + \ei \otimes \ei \in B \otimes B\,$. It is easy to see
that $\, \Delta_{B,S}(x) = 0 \,$ for all $ x \in S $, so $\,
\Delta_{B,S} \in \Der_S(B, \BB)\,$. This also follows immediately
from the fact that $ S $ is a separable algebra, and $\, \be \in S
\otimes S \,$ is the canonical separability element in $ S $.
\begin{lemma}
\la{DP} For any $\, \blambda \in S $, there is a canonical algebra
isomorphism
$$
\Pi^\blambda(B) \cong  T_{B}\Der_S(B, \BB)/ \langle \Delta_{B,S} -
\blambda \rangle\ .
$$
\end{lemma}

Thus, the structure of the algebras $\, \Pi^\blambda(B)\,$ is
determined by the bimodule $\,\Der_S(B, \BB)\,$. We now describe
this bimodule explicitly, in terms of $ \O \,$, $\, \I $ and the
dual line bundle $\, \I^* := \Hom_{\O}(\I,\,\O)\,$. To fix
notation we begin with a few fairly obvious remarks on bimodules
over one-point extensions.

A bimodule $ \Lambda $ over $ B = \O[\I\,] $ is characterized by
the following data: an $\O$-bimodule $ T $, a left $\O$-module $ U
$, a right $\O$-module $ V $ and a $\c$-vector space $ W $ given
together with three $\O$-module homomorphisms $\,f_1:\, \I \otimes
V \to T \,$, $\, f_2: \, \I \otimes W \to U \,$, $\, g_1:\, T
\otimes_\O \I \to U \,$ and a $\c$-linear map $\, g_2:\, V
\otimes_\O \I \to W \,$, which  fit into the commutative diagram
\begin{equation}
\la{dibim}
\begin{diagram}[small, tight]
\I \otimes  V \otimes_\O \I  & \rTo^{\id \otimes g_2}     & \I \otimes W \\
\dTo^{f_1 \otimes_\O \id}     &                            & \dTo_{f_2} \\
T \otimes_\O \I                & \rTo^{g_1}                 & U \\
\end{diagram}
\end{equation}
These data can be conveniently organized by using the matrix
notation
$$
\Lambda = \left(
\begin{array}{cc}
 T & U \\
 V & W
\end{array}
\right) \ ,
$$
with understanding that $ B $ acts on $ \Lambda $  by the usual
matrix multiplication, via the maps $\,f_1, \,f_2, \,g_1 $ and $
g_2\,$. Note that the commutativity of \eqref{dibim} ensures the
associativity of the action of $B$.

With this notation, the bimodule $\, \Der_S(B, \BB)\,$ can be
described as follows.
\begin{lemma}
\la{DerS}
There is an isomorphism of $B$-bimodules
$$
\Der_S(B, \BB) \cong
\left(
\begin{array}{cc}
 \Der(\O, \O^{\,\otimes 2})  & \Der(\O, \I\otimes \O)\\*[1ex]
 0 & 0
\end{array}
\right) \bigoplus \left(
\begin{array}{cc}
 \I \otimes \I^*  &  \I \otimes \O\\*[1ex]
 \I^* & \O
\end{array}
\right) \ ,
$$
with $\,\Delta_{B,S}\,$ corresponding to the element
$$
 \left[ \left(
\begin{array}{cc}
 \Delta  & 0\\*[1ex]
 0 & 0
\end{array}
\right) \, , \, \left(
\begin{array}{cc}
 -\sum_i \v_i \otimes \w_i  & 0\\*[1ex]
 0 & 1
\end{array}
\right) \right]\ ,
$$
where $\, (\v_i,\,\w_i) \,$ is a pair of dual bases\footnote{We
recall that $\, \{\v_i\} \subset \I \,$ and $\, \{\w_i\} \subset
\I^* \,$ form a `dual basis' for f.~g. projective modules $ \I $
and $ \I^*\,$ if $\,s = \sum_{i} \w_i(s)\,\v_i \,$ for all $\,s\in
\I \,$. Of course, this is an abuse of terminology, since $
\{\v_i\} $ and $\,\{\w_i\} \,$ are only generating sets of $ \I $
and $ \I^*$, not necessarily bases. The existence of such
generating sets characterizes f.~g. projective modules, see, e.g.,
\cite{B}, Ch.~II, Prop.~(4.5).} for the line bundles $ \I $ and $
\I^* $.
\end{lemma}

Now, as a consequence of Lemma~\ref{DP} and Lemma~\ref{DerS}, we get
\begin{proposition}[\cite{BC}, Proposition~5.1]
\la{gene}
The algebra $\,\Pi^\blambda(B)\,$ is generated by (the images of)
the following elements
$$
\ha := \left(
\begin{array}{cc}
 a & 0\\*[1ex]
 0 & 0
\end{array}
\right)\ ,\ \hv_i := \left(
\begin{array}{cc}
 0 & \v_i\\*[1ex]
 0 & 0
\end{array}
\right) \ , \ \hd :=\left(
\begin{array}{cc}
 d & 0\\*[1ex]
 0 & 0
\end{array}
\right)\ ,\ \hw_i := \left(
\begin{array}{cc}
 0 & 0\\*[1ex]
 \w_i & 0
\end{array}
\right) \ ,
$$
where $\, \ha \, ,\,\hv_i \in B \,$ and $\, \hd \,,\, \hw_i \in
\Der_S(B, \BB)\,$ with $\,d \in \Der(\O,\AA)\,$. Apart from the
obvious relations induced by matrix multiplication, these elements
satisfy
\begin{equation}
\la{r}
\hdel - \sum_{i=1}^N \hv_i \cdot \hw_i = \lambda\, e\
,\quad \sum_{i=1}^N \hw_i \cdot \hv_i = \lambdai\, \ei\ ,
\end{equation}
where `$\, \cdot $' denotes the action of $ B $ on the bimodule $
\Der_S(B, \BB)$.
\end{proposition}
With Proposition~\ref{gene}, we can describe the variety $
\CC_{n}(X, \I) $ as the space of equivalence classes of linear
maps (matrices)
$$
\{\,(\hha,\, \hhd,\, \hhv_i, \,\hhw_i) \ :\ \hha,\ \hhd \in \End(\c^n)
\, ,
\ \hhv_i \in
\Hom(\c,\c^n) , \  \hhw_i \in \Hom(\c^n, \c)\,\} \ ,
$$
satisfying the relations (cf. \eqref{E26})
\begin{equation}
\la{r0}
\bar{\Delta} - \sum_{i=1}^N \hhv_i \, \hhw_i = \id_n\
,\qquad \sum_{i=1}^N \hhw_i \, \hhv_i = - n \ .
\end{equation}
Of course, in addition to \eqref{r0}, $\,\hha\,$ and $\,\hhd \,$
should also obey the internal relations of the algebra $ \O $ and
the bimodule $ \Der(\O, \AA) $. Giving a matrix presentation of $
\CC_{n}(X, \I) $ thus boils down to describing $ \O $ and $
\Der(\O, \AA) $ in terms of generators and relations. This can be
easily done in practice.
\subsection{Example: plane curves}
\la{plane} Let $ X $ be a smooth curve on $ \c^2 $ defined by the
equation $\, F(x,y)=0\,$, with $\, F(x,y) :=
\sum_{r,s}a_{rs}x^ry^s \in \c[x, y]\,$. In this case, the algebra
$\, \O \cong \c[x,y] / \langle F(x,y)\rangle\,$ is generated by
$x$ and $y$, and the $\O$-module $\,\Der(\O)\,$ is (freely)
generated by the derivation $\,\partial\,$ defined by
\begin{equation*}
\partial(x)= F'_y(x,y)\,,\qquad \partial(y)= - F'_x(x,y)\ .
\end{equation*}

Further, it is easy to show that the bimodule $\,\Der(\O, \AA)\,$
is generated by the distinguished derivation $ \Delta =
\Delta_{\O} $ and the element $\,z\,$ defined by
\begin{gather*}
    z(x)=\sum_{r,s} a_{rs}\sum_{k=0}^{s-1}x^ry^k\otimes
    y^{s-k-1}\ ,\\
    z(y)=-\sum_{r,s} a_{rs}\sum_{l=0}^{r-1}x^l\otimes
    x^{r-l-1}y^s\ .
\end{gather*}
These generators satisfy the following commutation relations
\begin{gather}
\la{r1}
    [z,\, x]=\sum_{r,s} a_{rs}\sum_{k=0}^{s-1}y^{s-k-1}\Delta y^kx^r
    \,,\\\la{r2}
    [z,\, y]=-\sum_{r,s} a_{rs}\sum_{l=0}^{r-1}y^sx^{r-l-1}\Delta x^l
    \,.
\end{gather}
By Proposition~\ref{gene}, the algebra $\,\Pi^\blambda(B) \,$ is then
generated by the elements $\,\hat{x} $, $\,\hat{y} $, $\,\hat{z} $,
$\,\hv_i $, $\,\hw_i $ and $\, \hdel $,
subject to the relations \eqref{r} and \eqref{r1}, \eqref{r2}.

Let us now explicitly describe the generic points of the varieties
$ \CC_n(X, \I) $ (see remarks following Theorem~\ref{T6}). First,
we consider the case when $ \I $ is trivial, i.~e. $\,\I \cong \O
\,$. We choose $\,n\,$ distinct points $\, p_i = (x_i,\,y_i)\in
X\,$, $\,i=1,\dots, n\,$ and define the matrices
\begin{equation}
\la{ldata}
(\sX,\sY,\sZ,\sv,\sw) \in \End(\c^n) \times \End(\c^n)
\times \End(\c^n)\times \Hom(\c, \c^n)\times \Hom(\c^n,\c)
\end{equation}
by $\, \sX =\mathrm{diag}(x_1,\dots,x_n)\,$, $\ \sY
=\mathrm{diag}(y_1,\dots,y_n)\,$, $\ \sv^t = - \sw =
(1,\,\dots\,,\,1)\,$, and
\begin{equation}\la{mz}
    \sZ_{ii}=\alpha_i\quad ,\quad
    \sZ_{ij}= \frac{F(x_i,\, y_j)}{(x_i-x_j)(y_i-y_j)}
    \quad (i\ne j)\ ,
\end{equation}
where $\,\alpha_1,\dots ,\alpha_n \,$ are arbitrary scalars. A
straightforward calculation, using the relations \eqref{r1} and
\eqref{r2}, shows then that
\begin{equation*}
\hat{x} \mapsto \sX\ ,\quad  \hat{y} \mapsto \sY\ ,\quad \hat{z}
\mapsto \sZ\ ,\quad \hv \mapsto \sv\ , \quad \hw \mapsto \sw\ ,
\quad \hdel \mapsto \id_n + \sv\,\sw
\end{equation*}
defines a representation of $\,\Pi^\blambda(B)\,$ on the vector
space $\,\bV = \c^n \oplus \c\,$. The equivalence classes of such
representations correspond to generic points of
$\,\CC_n(X,\,\O)\,$.

\vspace*{1ex}

\remark\ The matrix $ \sZ $ defined by \eqref{mz} is a
generalization of the classical {\it Moser matrix} in the theory
of integrable systems (see \cite{KKS}).

\vspace*{1ex}

Now, let $ \I $ be an arbitrary line bundle on $ X $. As before,
we identify $\,\I\,$ with an ideal in $ \O $ and assume that the
zero set $ \V(\I) $ of $\I$ does not include the points $\, p_i
\,$ (this can always be achieved by changing the embedding of $ \I
$ in $ \O $). Then, $ \I^* $ can be identified with a fractional
ideal of $ \O $ generated by rational functions with poles in $
\V(\I) $, and the pairing $\, \I \times \I^* \to \O\,$ is given by
multiplication in $ \c(X) $. The evaluation of $\, \v \in \I \,$
at $\,p_1,\,\dots,\, p_n\,$ defines a vector $\,\hhv \in \c^n $;
in a similar fashion, any $\,\w \in \I^* $ defines a row vector $
\hhw =(\hhw_1\,,\dots , \hhw_n)\in \Hom(\c^n,\, \c) $, with
$\hhw_j=-\w(p_j)$. If $\{\v_i\}$, $\{\w_i\}$ are dual bases for $
\I $ and $ \I^* $, then  $\, \sum_{i} \v_i \otimes \w_i \,$ gives
a rational function on $ X \times X $, which we denote by $ \phi
$; the fact that the bases are dual implies $\,\phi(p,\,p)=1\,$
for all $ p \in X $. As a result, the $n\times n$ matrix
$\,\sum_{i}\hhv_i\, \hhw_i \,$ equals $ \,\|-\phi(p_i,\,
p_j)\|\,$, with all the diagonal entries being equal to $-1$.

Now, let $ \sX $ and $\, \sY $ be the diagonal matrices as above,
and let $ \sZ $ be given by
$$
    \sZ_{ii}=\alpha_i\quad ,\quad
    \sZ_{ij}= \frac{F(x_j,\, y_i)\,\phi(p_i,\, p_j)}{(x_i-x_j)(y_i-y_j)}
    \quad (i\ne j)\ .
$$
It is straightforward to check that the assignment
\begin{equation*}
\hat{x} \mapsto \sX\ ,\quad  \hat{y} \mapsto \sY\ ,\quad \hat{z}
\mapsto \sZ\ ,\quad \hv_i \mapsto \sv_i\ , \quad \hw_i \mapsto
\sw_i\ , \quad \hdel \mapsto \id_n + \sum_i \sv_i\,\sw_i
\end{equation*}
extends to a representation of $\Pi^\blambda(B)$ on the vector
space $ \bV = \c^n \oplus \c \,$; such representations correspond
to generic points of the variety $ \CC_n(X,\I) $.

\section{The Calogero-Moser Correspondence}
\la{CMMap} We will keep the notation from the previous section: $
\O = \O(X) $ stands for the coordinate ring of a smooth affine
irreducible curve $X$,  $\, B = \O[\I]\,$ for the one-point
extension of $ \O $ by a line bundle $\, \I \,$, and $\, \Pi =
\Pi^\blambda(B) \,$ for the deformed preprojective algebra over $
B $ of weight $ \blambda = (1,-n) $ with $ n \in \N $.

\subsection{Recollement}
We now explain the relation between representations of $
\Pi^\blambda(B) $ and $\D$-modules on $X$. We begin with the
following observation, which is a simple consequence of
Theorem~\ref{TCB1} (see \cite{BC}, Lemma~4.1).
\blemma \la{recoll} There is a canonical algebra map $\, \btheta:
\Pi^{\blambda}(B) \to \Pi^{1}(\O) \,$, which is a surjective
pseudo-flat ring epimorphism, with $\, \Ker(\btheta) = \langle\,
\ei\,\rangle\,$. \elemma

With Proposition~\ref{gene}, the homomorphism $\, \btheta \,$  can
be described explicitly, in terms of generators of $
\Pi^\blambda(B)\, $:\, namely,
\begin{equation}
\la{bth}
\btheta(\ha) = \overline{a}\ , \quad \btheta(\hd) = \overline{d} \ ,
\quad \btheta(\hv_i) = \btheta(\hw_i) = 0\ ,
\end{equation}
where $\, \overline{a} \,$ and $\, \overline{d} \,$ denote the
classes of $\, a \in \O \,$ and $\, d \in \Der(\O, \AA) \,$ in the
tensor algebra of $\, \Der(\O, \AA) \,$ modulo the ideal $\,
\langle \Delta - 1 \rangle \,$.

Now, by Example~\ref{TCB2} (see Sect.~\ref{DPA}), the algebra $\,
\Pi^{1}(\O)\,$ is isomorphic to $\,\D\,$: we fix, once and for
all, such an isomorphism to identify $ \D = \Pi^{1}(\O) $. In
combination with Lemma~\ref{recoll}, this yields an algebra map
$\,\btheta:\,\Pi \to \D \,$. We will use $ \btheta $ to relate the
{\it derived} module categories of $ \Pi $ and $ \D $ in the
following way (cf. \cite{BCE}).

First, we let $\, U \,$ denote the endomorphism ring of the
projective module $\, \ei \Pi \,$: this ring can be identified
with the associative subalgebra $\, \ei \Pi\, \ei \,$ of $\,\Pi
\,$ having $ \ei $ as an identity element. Next, we introduce six
additive functors $ (\btheta^*, \,\btheta_*,\, \btheta^{!}) $ and
$ (j_!, \,j^*,\, j_*) $ between the module categories of $ \Pi $,
$\,\D$ and $ U $. We define $\, \btheta_*:\,\Mod(\D)\to \Mod(\Pi)
\,$ to be the restriction functor associated to the algebra map
$\, \btheta:\,\Pi \to \D \,$. This functor is fully faithful and
has both the right adjoint $\, \btheta^! := \Hom_{\Pi}(\D,\,
\mbox{---}\,)\,$ and the left adjoint $\, \btheta^* := \,\D
\otimes_{\Pi} \mbox{---}\,$, with adjunction maps $\, \btheta^*
\btheta_* \simeq \id \simeq \btheta^!\,\btheta_* \,$ being
isomorphisms. Now we define $\, j^*:\, \Mod(\Pi) \to \Mod(U) \,$
by $\, j^*\bV := \ei \bV \,$. Since $ e_\infty \in \Pi $ is an
idempotent, $ j^* $ is exact and has also the right and the left
adjoint functors: $\, j_* := \Hom_{U}(e_\infty \Pi,
\,\mbox{---}\,)\,$ and  $\, j_! := \Pi \ei \otimes_{U}
\mbox{---}\,$ respectively, satisfying $\, j^*j_* \simeq \id
\simeq j^* j_! \,$.

The functors $ (\btheta^*, \,\btheta_*,\, \btheta^{!}) $ and $
(j_!, \,j^*,\, j_*) $ defined above induce the six exact functors
at the level of (bounded) derived categories:
\begin{equation}
\la{Dia1}
\begin{diagram}[small, tight]
\DC^{-}(\Mod\, \D) \ & \ \pile{\lTo^{L \btheta^*}\\
\rTo^{\btheta_{\!*}} \\
\lTo^{\RR \btheta^!}} \ &
\ \DC^{-}(\Mod\, \Pi)\ & \ \pile{\lTo^{L j_!}\\
\rTo^{j^*}\\ \lTo^{\RR j_*}} \ & \ \DC^{-}(\Mod\, U) \\
\end{diagram}\ .
\end{equation}
The properties of these functors can be summarized in the
following way.
\bprop[\cite{BC}, Proposition~4.1] \la{P1} The diagram
\eqref{Dia1} is a recollement of triangulated categories.
\eprop

\remark\ The `recollement axioms' were originally formulated in
\cite{BBD} to imitate a natural structure on the derived category
$ \DC({\mathscr Sh}_X) $ of abelian sheaves arising from the
stratification of a topological space into a closed subspace and
its open complement (see \cite{BBD}, Sect.~1.4). In an algebraic
setting similar to ours, these axioms were studied in \cite{CPS}.

We will use the induction functor $\, L\btheta^* \,$ to define a
natural map: $\,\CC_n(X, \I) \to \R(\D) \,$. As a first step, we
compute the values of $\, L\btheta^* \,$ on the finite-dimensional
representations of $ \Pi^\blambda(B)$, regarding the latter
as $0$-complexes in $\DC^-(\Mod\,\Pi)\,$ (see \cite{BC},
Lemma~4.2). We recall that $\,L_n \btheta^* \cong
\Tor^{\Pi}_n(\D,\,\mbox{---}\,) \,$, where $ \D $ is viewed as a $
\Pi$-module via the algebra map $ \btheta $.

\blemma \la{simp} If $ \bV $ is a $\Pi$-module of finite dimension
over $\,\c\,$, then $\, L_n \btheta^* (\bV) = 0 \,$ for $\, n
\not= 1 \,$ and
\begin{equation}
\la{CMf}
L_{1} \btheta^* (\bV) \cong \Ker\left[\,\Pi\,
\ei \otimes_{U} \ei \bV \xrightarrow{\bmu} \bV \right],
\end{equation}
where $ \bmu $ is the natural multiplication-action map on $\bV$.
\elemma
Lemma~\ref{simp} shows that $\,L\btheta^*(\bV)\,$ is isomorphic in
$\, \DC^{-}(\Mod\, \D)\,$ to a single $\D$-module $\,M = L_{1}
\btheta^* (\bV) \,$ located in cohomological degree $(-1)$.
Abusing notation, we will simply write $\, M = L\btheta^*(\bV) \,$
in this case.

\subsection{The action of $ \Pic(\D) $ on Calogero-Moser spaces}
\la{ad}
Next, we describe a natural action of the Picard group of
the category of $\D$-modules on representation varieties of $
\Pi^\blambda(B) $. It is known that $\, \Pic(\D) \,$ has different
descriptions depending on whether $ X $ is the affine line or not
(see \cite{CH1}). In this section, we will assume that $ X \ne
\A^1 $. Our main result (Theorem~\ref{Tmain}) will still be true
for all $X$, since the case $ X = \A^1 $  is covered in \cite{BW1,
BW2}.

We recall (see, e.g., \cite{B}, Ch.~II, \S~5) that $ \Pic(\D) $
can be identified with the group of $ \c$-linear auto-equivalences
of the category $ \Mod(\D)$, and thus it acts naturally on $
\R(\D) $ and $ K_0(\D) $. To be precise, the elements of $\,
\Pic(\D) \,$ are the isomorphism classes $ [\P] $ of invertible
bimodules over $\D$, and the action of $ \Pic(\D) $ on $ \R(\D) $
and $ K_0(\D) $ is defined by $\, [M] \mapsto [\P \otimes_{\D}
M]\,$. Observe that $ \Pic(\D) $ acts on $ K_0(\D) $ preserving
rank: hence, this action restricts to $ \Pic(X) $ through the
identification $\,K_0(\D) \cong K_0(X) \cong \Z \oplus \Pic(X)
\,$, see Section~\ref{SCI}.

Now, let $ \Aut(\D) $ denote the group of $\c$-automorphisms of
the algebra $ \D $. There is a natural homomorphism: $\,\Aut(\D)
\to \Pic(\D) \,$, sending $\,\vp \in \Aut(\D) \,$ to (the class
of) the bimodule $ \D_\vp $. The kernel of this homomorphism
consists of the inner automorphisms of $ \D$, so the group
$\,\Gamma := \Out(\D)\,$ of {\it outer} automorphism classes can
be identified with a subgroup of $ \Pic(\D) $. With this
identification, we have
\bprop[see \cite{BW}, Theorem~1.1] \la{piceq} $\, \Pic(\D) $ acts
on $ \Pic(X) $ transitively, the stabilizer of a point $\,[\I] \in
\Pic(X) $ being isomorphic to $\,\Gamma\,$. The map $
\gamma:\,\R(\D) \to \Pic(X) $ defined by \eqref{E10} is
equivariant under the action of $ \Pic(\D)$. \eprop

Explicitly, we can describe the action of $ \Pic(\D) $ on $
\Pic(X) $ as follows. By \cite{CH1}, every invertible bimodule
over $ \D $ is isomorphic to $\, \D\L = \D \otimes_{\O} \L \,$ as
a left module, while the right action of $ \D $ on $ \D\L $ is
determined by an algebra isomorphism $\, \vp:\, \D
\stackrel{\sim}{\to} \End_{\D}(\D\L) \,$, where $\, \L \,$ is a
line bundle on $X$. Following  \cite{BW}, we write $ \,
(\D\L)_{\vp} \,$ for this bimodule. Restricting $ \vp $ to $\O$
yields an automorphism of $X$, and the assignment
\begin{equation}
\la{Egr}
g:\,\Pic(\D) \to \Pic(X) \rtimes \Aut(X)\ ,\quad
[\, (\D\L)_{\vp}] \mapsto ([\L],\, \vp |_\O)\ ,
\end{equation}
defines then a group homomorphism. On the other hand, $\, \Pic(X)
\rtimes \Aut(X)\,$ acts on $ \Pic(X) $ in the obvious way, via
left multiplication:
\begin{equation}
\la{E13}
([\L],\tau) : \, [\,\I\,] \mapsto [\,\L \,\tau(\I)\,] \ ,
\end{equation}
where $\, ([\L],\tau) \in  \Pic(X) \rtimes \Aut(X)\,$ and $\, [\I]
\in \Pic(X) \,$. Combining \eqref{Egr} and \eqref{E13} together,
we get an action of $ \Pic(\D) $ on $ \Pic(X) $, which is easily
seen to agree with the natural action of $ \Pic(\D)$ on $
K_0(\D)$.

Now, given a line bundle $ \I $ and an invertible bimodule $\P =
(\D\L)_{\vp}\, $, we set
$$
P := \tilde{\L} \otimes_{\tA} B_{\tau} \ ,
$$
where $\,\tA := \O \times \c \,$, $\, \tilde{\L} := \L \times \c
\,$ and $\, B_{\tau} := \O[\tau(\I)] \,$ with $\, \tau =
\vp|_{\O}\,$. By Lemma~\ref{L10}$(a)$, $\, P $ is a progenerator
in the category of right $ B_\tau$-modules, with endomorphism ring
$$
\End_{B_\tau}(P) \cong \tilde{\L} \otimes_{\tA} B_\tau
\otimes_{\tA} \tilde{\L}^* \cong \O[\J] \ ,
$$
where  $\, \tilde{\L}^* := \L^* \times \c \,$ and $\, \J :=
\L\,\tau(\I) \,$. Thus, associated to the bimodule $ \P $ is the
Morita equivalence:
$$
\Mod(B_\tau) \stackrel{\sim}{\to} \Mod(\O[\J])\ , \quad
\bV \mapsto P \otimes_{B_\tau} \bV \cong \tilde{\L} \otimes_{\tA}
\bV\ .
$$
Next, we extend $\, P \,$ to a $ \Pi^\blambda(B_\tau)$-module by
\begin{equation}
\la{prog}
\bP := P \otimes_{B_\tau} \Pi^\blambda(B_\tau) \cong
\tilde{\L} \otimes_{\tA} \Pi^\blambda(B_\tau)\ ,
\end{equation}
which is clearly a progenerator in the category of right $
\Pi^\blambda(B_\tau)$-modules. By Lemma~\ref{L10}$(b)$, the
algebra $ B_{\tau} $ is isomorphic to $ B \,$: the isomorphism is
given by
\begin{equation}
\la{tiso}
\ttau:\ B \to B_\tau\ ,\quad
\left(
\begin{array}{cc}
 a  & b\\
 0  & c
\end{array}
\right) \mapsto
\left(
\begin{array}{cc}
 \tau(a) & \tau(b) \\
 0  & c
\end{array}
\right)\ .
\end{equation}
Since $\,\ttau(\blambda) = \blambda \,$ for all $\,\blambda \in S
\,$, \, \eqref{tiso} canonically extends to an isomorphism of
deformed preprojective algebras: $\,\Pi^\blambda(B)
\stackrel{\sim}{\to} \Pi^\blambda(B_\tau) \,$, which we will also
denote by $ \ttau $. Now, using this last isomorphism, we will
regard $ \bP $ as a $ \Pi^\blambda(B)$-module and identify
\begin{equation}
\la{endp}
\End_{\Pi^\blambda(B)}(\bP) =
\tilde{\L} \otimes_{\tA} \Pi^\blambda(B_\tau)
\otimes_{\tA}\tilde{\L}^* \cong
\tilde{\F} \otimes_{\tA} \Pi^\blambda(B)
\otimes_{\tA}\tilde{\F}^*\  ,
\end{equation}
where $\,\F := \L^\tau = \tau^{-1}(\L)\,$ and $ \tilde{\F} = \F
\times \c $. With identification \eqref{endp}, we have the
embedding
\begin{equation}
\la{lso}
\ttau^{-1}:\ A[\J] \cong \tilde{\L} \otimes_{\tA} B_\tau
\otimes_{\tA} \tilde{\L}^* \into \End_{\Pi^\blambda(B)}(\bP)\ ,
\end{equation}
and, since $\,\End_{\D}(\F \D) = \F \otimes_A \D \otimes_A \F^*
\cong \tilde{\F} \otimes_{\tA} \D \otimes_{\tA} \tilde{\F}^* \,$,
the natural map
\begin{equation}
\la{lso1}
1 \otimes \btheta \otimes 1 :\ \End_{\Pi^\blambda(B)}(\bP) \to
\End_{\D}(\F \D)\ ,
\end{equation}
where $\, \btheta:\,\Pi^\blambda(B) \to \D \,$ is given by
Lemma~\ref{recoll}.

On the other hand, $\, \vp(\D) = \End_{\D}(\D\L) = \L^*\, \D\,  \L
\, $ implies $\,\D = \L\,\vp(\D)\,\L^*\,$, so taking the inverse
of $\, \vp \,$ defines an algebra isomorphism
\begin{equation}
\la{lso2}
\psi = \vp^{-1}:\ \D \to \F \,\D \,\F^* =
\End_{\D}(\F \D)\ .
\end{equation}
Combining \eqref{lso1} and \eqref{lso2} together, we get the
diagram of algebra maps
\begin{equation}
\la{pil}
\begin{diagram}
 \Pi^\blambda(\O[\J])   &  \rDotsto{ \ } & \End_{\Pi^\blambda(B)}(\bP) \\
\dTo^{\btheta}          &       & \dTo_{1 \otimes \btheta \otimes 1} \\
\D       & \rTo^{\psi\quad}       & \End_{\D}(\F \D) \\
\end{diagram}
\end{equation}
which obviously commutes when the dotted arrow is restricted to
\eqref{lso}.

\bprop[\cite{BC}, Proposition~4.3] \la{Lex} There is a {\rm
unique} algebra isomorphism $\, \bvp :\,\Pi^\blambda(\O[\J]) \to
\End_{\Pi^\blambda(B)}(\bP)\,$, extending \eqref{lso} and making
\eqref{pil} a commutative diagram. \eprop

The proof of this result in \cite{BC} is rather indirect: it
combines homological arguments of Theorem~\ref{lift} with explicit
calculations and description of automorphisms of $ \D $ given in
\cite{CH1}.

Now, using the isomorphism $ \bvp $ of Proposition~\ref{Lex}, we
can make $ \bP $ a left $ \Pi^\blambda(A[\J])$-module and thus a
progenerator from $ \Pi^\blambda(A[\I])$ to $ \Pi^\blambda(A[\J])
$. This assigns to $ \P = (\D\L)_{\vp} $ the Morita equivalence:
$$
\Mod\,\Pi^\blambda(A[\I]) \to \Mod\,\Pi^\blambda(A[\J])\ ,\quad
\bV \mapsto \bP \otimes_{\Pi} \bV\ ,
$$
which, in turn, induces an isomorphism of representation varieties
\begin{equation}
\la{eigg}
f_{\P}:\ \CC_n(X, \I) \stackrel{\sim}{\to} \CC_n(X, \J)\ .
\end{equation}

We warn the reader that \eqref{eigg} {\it depends} on the choice
of a specific representative in the class $\, [\P] \in \Pic(\D)
\,$, so, in general, we do not get an action of $ \Pic(\D) $ on
$\,\bigsqcup_{[\I] \in \Pic(X)}\,\CC_n(X,\I)\,$. However, it turns
out that $ f_{\P} $ induce a well-defined action of $ \Pic(\D) $
on the {\it reduced} Calogero-Moser spaces
$\overline{\CC}_n(X,\I)$. Precisely, we have

\blemma[\cite{BC}, Lemma~4.3] \la{omeac} The map \eqref{eigg}
induces an isomorphism
\begin{equation*}
\la{eiq}
\bar{f}_{\P}:\ \overline{\CC}_{n}(X,\I) \stackrel{\sim}{\to}
\overline{\CC}_{n}(X,\J)\ ,
\end{equation*}
which depends only on the class of $ \,\P $ in $ \Pic(\D) $.
\elemma
If we set $\, \overline{\CC}_n(X) := \bigsqcup_{[\I]\in
\Pic(X)}\,\overline{\CC}_{n}(X,\I)\,$, the assignment $\,[\P]
\mapsto \bar{f}_{\P} \,$ defines an action of $ \Pic(\D) $ on $
\overline{\CC}_{n}(X)\,$  for each $\,n \ge 0\,$.

\subsection{The action of automorphisms}
\la{actaut} Assume now that $\, \P \,$ comes from an automorphism
of $ \D $, i.~e. $\,[\P] \in \Gamma \subseteq \Pic(\D)\,$, where
$\,\Gamma := \Out(\D)\,$. Then, by Proposition~\ref{piceq}, $\,
[\P] \,$ stabilizes $\, [\I]\in \Pic(X)\,$, so the isomorphisms $
\bar{f}_{\P} $ of Lemma~\ref{omeac} define an action of $ \Gamma $
on each space $ \overline{\CC}_n(X, \I) $ individually. We now
describe this action in explicit terms.

First of all, when $\, X \ne \A^1$, we can identify (see
\cite{BC}, Section~5.5)
$$
\Gamma \cong \Aut(X) \ltimes (\Omega^1 X)/\Lambda \ ,
$$
where $ \Omega^1 X $ is the canonical bundle of $ X $ and $\,
\Lambda := \O^{\times}\!/\c^{\times} \,$ is the multiplicative
group of (nontrivial) units of $ \O(X) $ embedded in $\,\Omega^1
X\,$ via the logarithmic derivative map:
\begin{equation}
\la{uuu}
{\tt dlog}:\ \Lambda \into \Omega^1 X\,
\quad u \mapsto u^{-1} du\ .
\end{equation}
To simplify calculations, we will assume here that $ \Aut(X) $ is
trivial, which is clearly the case for generic curves.

Let $ \Omega^1(B) $ denote the bimodule of noncommutative
differential forms on $ B $ (i.e. the kernel of the multiplication
map $\,m:\,B \otimes B \to B \,$), and let $\, \DB :=
\Omega^1(B)/[B, \,\Omega^1 B]\,$ be the quotient of $ \Omega^1(B)
$ by its commutator subspace (the Karoubi-de Rham differentials of
$ B $). Using the fact that $ B $ is finitely generated and
quasi-free  (see Lemma~\ref{LL6}), we identify
\begin{equation}
\la{Kar}
\DB \cong B \otimes_{\eB} \Omega^1(B) \cong B
\otimes_{\eB} (\Omega^1 B)^{\star\star} \cong \Hom_{\eB}((\Omega^1
B)^{\star},\,B)\ ,
\end{equation}
where $\,\eB := B \otimes B^{\mbox{\tiny opp}}\,$, and $\,
(\mbox{---} \,)^\star \,$ stands for the duality over $ \eB $.
Explicitly, under the identification \eqref{Kar}, $\,
\overline{\omega} \in \DB \,$ corresponds to the map $\,
\hat{\omega}:\,\delta \mapsto m[\delta(\omega)]\,$.

Now, we define an action of $\, \DB \,$ on  $\, T_B\, (\Omega^1
B)^{\star} $ as follows: if $\, \overline{\omega} \in \DB \,$, we
let $\, \tsigma_\omega \,$ denote the automorphism of $\, T_B\,
(\Omega^1 B)^{\star}\,$ acting identically on $ B $ and mapping
$$
(\Omega^1 B)^{\star} \to B
\oplus (\Omega^1 B)^{\star} \into T_B\, (\Omega^1 B)^{\star}\ ,\quad
\delta \mapsto \delta + \hat{\omega}(\delta) \ .
$$
By the universal property of tensor algebras, this uniquely
determines $ \tsigma_\omega $, and it is clear that this map is
bijective. Moreover, if $ \Delta_B \in (\Omega^1 B)^\star $ is the
canonical inclusion $ \Omega^1 B \into \eB $, then $
\hat{\omega}(\Delta_B) = 0 $, and hence $ \tsigma_\omega(\Delta_B)
= \Delta_B $ for any $ \overline{\omega} \in \DB $. The assignment
$\, \overline{\omega} \mapsto \tsigma_\omega \,$ defines thus a
homomorphism
\begin{equation}
\la{auto1} \tsigma: \DB \to \Aut_{B}[T_B\, (\Omega^1 B)^{\star}]
\end{equation}
from the additive group of $ \DB $ to the subgroup of $B$-linear
automorphisms of $ T_B (\Omega^1 B)^{\star} $ preserving the
element $ \Delta_B $.

Next, identifying the canonical bundle of $ X $ with the module of
K\"ahler differentials of $ \O$, we construct an embedding of $
\Omega^1 X $ into $ \DB $. For this, we consider the exact
sequence
\begin{equation}
\la{auto2} 0 \to \rm{HH}_1(B) \stackrel{\alpha}{\to} \DB \to B \to
\rm{HH}_0(B) \to 0\ ,
\end{equation}
obtained by tensoring the fundamental exact sequence
$$
 0 \to \Omega^1(B) \to \eB \to B \to 0
$$
with $ B $, and compose the connecting map $ \alpha $ in
\eqref{auto2} with natural isomorphisms
\begin{equation}
\la{auto3} \rm{HH}_1(B) \cong \rm{HH}_1(\O) \cong \Omega^1 X \ .
\end{equation}
(The first isomorphism in \eqref{auto3} is induced by the
projection $\theta:\, B \to \O\,$, while the second by the
canonical map: $\, \AA \to \Omega^1 X \,$, $\,f \otimes g \mapsto
f \, dg $.)

Finally, combining \eqref{auto1} with \eqref{auto2} and
\eqref{auto3}, we define
\begin{equation}
\la{auto4}
\sigma:\ \Omega^1 X \stackrel{\alpha}{\into} \DB
\stackrel{\tsigma}{\to} \Aut_{B}[T_B\, (\Omega^1
B)^{\star}] \to \Aut_{S}[\Pi^\blambda(B)]\ ,
\end{equation}
where the last map is induced by the algebra projection: $\,T_B\,
(\Omega^1 B)^{\star} \onto \Pi^\blambda(B)\,$. With
Proposition~\ref{gene}, the action \eqref{auto4} can be easily
described in terms of generators of $ \Pi^\blambda(B) $.
\blemma \la{actgen} The homomorphism $\,\sigma:\,\Omega^1 X \to
\Aut_S[\Pi^\blambda(B)]\,$ is determined by
$$
\sigma_\omega(\ha) = \ha\ ,\quad \sigma_\omega(\hv_i) =
\hv_i \ ,\quad \sigma_\omega(\hw_i) = \hw_i \ ,\quad
\sigma_\omega(\hd) = \hd + \widehat{\omega(d)} \ ,
$$
where $\, \omega \in \Omega^1 X \,$ acts on $\, d \in \Der(\O,
\AA) \,$ via the natural identification, cf. \eqref{Kar}:
$$
\Omega^1 X \cong \DA \cong \Hom_{\O_{X \times X}}(\Der(\O, \,\AA),\,\O)
 \ .
$$
\elemma

Now, the group $\, \Aut_{S}[\Pi^\blambda(B)] \,$ acts on
$\,\Rep_S(\Pi^\blambda(B), \bn)\,$ in the obvious way: if
$\,\varrho: \Pi^\blambda(B) \to \End(\bV)\,$ represents a point in
$\,\Rep_S(\Pi^\blambda(B), \bn)\,$, then $\,(\omega, \varrho)
\mapsto \varrho \, \sigma_\omega^{-1} \,$ for $ \omega \in
\Omega^1 X $. Clearly, this commutes with the natural $
\GL_S(\bn)$-action on $\,\Rep_S(\Pi^\blambda(B), \bn)\,$ and hence
induces the action
\begin{equation}
\la{act1} \sigma^*:\ \Omega^1 X \to \Aut\, [\,\CC_{n}(X, \I)\,]
\ ,\quad \omega \mapsto [\,\sigma^*_\omega:\,\varrho \mapsto
\varrho \, \sigma_\omega^{-1}\,]
\ .
\end{equation}
By a straightforward calculation we get (see \cite{BC}, Lemma~4.3)
\blemma \la{ommm} The action \eqref{act1} agrees with
\eqref{eigg}{\rm\,:\,} that is, if $\,\P = \D_{\sigma_\omega}\,$,
then
$$
f_{\P} = \sigma^*_\omega \quad \mbox{for all} \quad
\omega \in \Omega^1 X\ .
$$
\elemma

On the other hand, the restriction of \eqref{act1} to the group
$\,\Lambda\,$ via \eqref{uuu} agrees with the natural action of
$\, \Aut_X(\I) = \O^\times \,$ on $\,\CC_{n}(X, \I)\,$. Thus
\eqref{act1} induces an action of $\,\Gamma = (\Omega^1
X)/\Lambda\,$ on each of the spaces $ \overline{\CC}_n(X, \I) $.
By Lemma~\ref{ommm}, this coincides with the restriction to $
\Gamma $ of the action of $ \Pic(\D) $ constructed in
Section~\ref{ad}.

\subsection{The main theorem}
\la{TMT} We can now state the main result of \cite{BC}. We recall
the functor $\, L\btheta^* = \Tor^\Pi_1(\D,\,\mbox{---}):\,
\Mod(\Pi) \to \Mod(\D) \,$ associated to the algebra homomorphism
$\,\btheta:\,\Pi \to \D \,$: when restricted to finite-dimensional
representations, this functor is given by \eqref{CMf}.
\bthm[\cite{BC}, Theorem~4.2] \la{Tmain} Let $X$ be a smooth
affine irreducible curve.

$(a)$ For each $ n \geq 0 $ and $ [\I] \in \Pic(X) $, the functor
$\,L\btheta^* $ induces an injective map
$$
\omega_n:\,\overline{\CC}_n(X, \I) \to \gamma^{-1}[\I]\ ,
$$
which is equivariant under the action of the group $\,\Gamma $.

$(b)$ Amalgamating  the maps $\, \omega_n \,$ for all $\, n \ge
0\, $ yields a bijective correspondence
$$
\omega:\ \bigsqcup_{n \geq 0} \overline{\CC}_n(X,\I)
\stackrel{\sim}{\to} \gamma^{-1}[\I] \ .
$$

$(c)$ For any $ [\I] $ and $ [\J] $ in $ \Pic(X) $ and for any
$[\P] \in \Pic(\D) $, such that $ [\P] \cdot [\I] = [\J] $, there
is a commutative diagram:
\begin{equation}
\la{CMdim}
\begin{diagram}[small, tight]
 \overline{\CC}_n(X,\I)    &  \rTo^{\ \bar{f}_{\P}\ } & \overline{\CC}_n(X,\J) \\
\dTo^{\omega_n}            &                   & \dTo_{\omega_n} \\
\gamma^{-1}[\I]            & \rTo^{\ [\P]\ }       & \gamma^{-1}[\J] \\
\end{diagram}
\end{equation}
where $ \bar{f}_{\P} $ is an isomorphism induced by \eqref{eigg}.

\ethm

For technical reasons, we assumed above that $\, X \ne \A^1 $.
However, Theorem~\ref{Tmain} holds true in general: if $ X = \A^1
$, by \cite{BCE}, Theorem~1, the map $ \omega $ induced by $
L\btheta^* $ agrees with the Calogero-Moser map constructed in
\cite{BW1, BW2}. In this case, the ring $ \D $ is isomorphic to
the Weyl algebra $ A_1(\c) $ and $\, \Pic(\D) \cong \Aut(\D) =
\Out(\D) $, see \cite{St}.  Since $ \Pic(\A^1) $ is trivial, the
last part of Theorem~\ref{Tmain} implies the equivariance of $
\omega $ under the action of $ \Aut(A_1) $, which was one of the
main results in \cite{BW1}.

\subsection{Explicit construction of ideals}
\la{weyl} To illustrate Theorem~\ref{Tmain}, we return to our
basic example of plane curves (see Section~\ref{plane}). In this
case, we will describe the map $ \omega $ quite explicitly, in
terms of generalized Calogero-Moser matrices \eqref{ldata}. For
simplicity, we consider only the case when $ \I $ is trivial. A
$\Pi$-module $\, \bV = \c^n \oplus \c \,$ can then be described by
$ (\sX,\,\sY,\,\sZ, \,\sv,\, \sw)$, which, apart from \eqref{r1}
and \eqref{r2}, satisfy the following relations
$$
F(\sX,\sY)=0\ ,\quad [\sX,\,\sY] = 0\quad \mbox{and}\quad \sdel
= \id_n+\sv\,\sw \ .
$$
The dual representation $\,\varrho^*:\, \Pi^{\rm \footnotesize
opp} \to \End(\bV^*) \,$ is given by the transposed matrices $
(\sX^t,\,\sY^t,\,\sZ^t, \,\sv^t,\, \sw^t)$.

Now, using these matrices, we define the following element in the
field of fractions of the algebra $ \D\,$:
\begin{equation*}
\la{kappz}
\kappa :=1+\sv^t\,(\sZ^t-z\,\id_n)^{-1}(\sX^t-x\,\id_n)^{-1}
(\sY^t-y\,\id_n)^{-1} F(\sX^t,\, y\,\id_n)\, \sw^t \ .
\end{equation*}
and consider the (fractional) left ideal
\begin{equation*}\la{pca}
M_{\bV} := \D \, \det(\sX-x\,\id_n) + \D \, \det(\sY-y\,\id_n) +
\D\, \det(\sZ-z\,\id_n)\,\kappa\ .
\end{equation*}
\bprop \la{exxp} If $\, [\bV] \in \CC_n(X, \I)\,$ is determined by
$\,(\sX,\,\sY,\sZ,\,\sv,\,\sw)\,$, the corresponding ideal class
$\,\omega[\bV] \,$ in $\,\R(\D)\,$ is represented by $ M_{\bV} $.
\eprop
For the proof of Proposition~\ref{exxp} and more examples, we
refer the reader to \cite{BC}, Section~6. Here, we only mention
that Theorem~6.1 of \cite{BC} gives a similar explicit
construction of $\, \omega \,$ for an {\it arbitrary} smooth
curve, generalizing earlier calculations for the first Weyl
algebra in \cite{BC1}.

\bibliographystyle{amsalpha}

\end{document}